\documentclass[11pt]{article}
\usepackage{amssymb}
\newtheorem{lemma}{Lemma}[section]

\newtheorem{theorem}[lemma]{Theorem}
\newtheorem{proposition}[lemma]{Proposition}
\newtheorem{corollary}[lemma]{Corollary}
\newtheorem{example}[lemma]{Example}

\newtheorem{remark}[lemma]{Remark}

\def\endproof{\hfill$\Box$}

\def\endproof{\hfill$\Box$}

\title{Kinds of preservations for properties\footnote{The work was carried out in the framework of the
State Contract of the Sobolev Institute of Mathematics, Project
No.~FWNF-2022-0012.}}
\author{T.E. Rajabov, S.V. Sudoplatov}
\date{}

\begin{document}

\maketitle
\begin{abstract}
We study possibilities of preservations for properties, their
links and related connections between semantic and syntactic ones,
both in general and as characterizations for subalgebras,
congruence relations, Henkin construction, Tarski-Vaught test,
some graphs, and graded structures. Traces for types and related
objects are studied, too.
\end{abstract}


\bigskip
\section{Introduction}

We study various properties \cite{rs1} which are preserved under
given conditions. These preservations generalize the notion of
$(p,q)$-preserving formula \cite{rs2, rs3, rs4} and its variations
for correspondent type-definable sets.

The paper is organized as follows. In Section 2, we introduce
various kinds of preservation and describe various related
properties. Traces of types are introduced in terms preservations
and their general properties are found. In Section 4, subalgebras
and congruence relations are characterized in terms of
preservations. In Section 5, connections of model constructions,
including Henkin construction, are described in terms of
preservations of formulae. Vaught-Tarski text is characterized in
terms of preservations in Section 6. Multipartite and related
graphs are connected with preservations in Section 7. In Section
8, we study preservations of properties in graphs with respect to
distances. Preservation properties for algebraic constructions and
decompositions are discussed in Section 9.

Throughout we use the standard model-theoretic notions and
notations.

\section{Preservations of properties by types and formulae, their general properties}

{\bf Definition.} Let $\mathcal{M}$ be a structure, $P_1\subseteq
M^{k_1},\ldots,P_n\subseteq M^{k_n}$, $Q\subseteq M^m$ be
properties,
$\Phi=\Phi(\overline{x}_1,\ldots,\overline{x}_n,\overline{y})$ be
a type with $l(\overline{x}_1)=k_1,$ $\ldots,$
$l(\overline{x}_n)=k_n$, $l(\overline{y})=m$. We say that the
tuple $(P_1,\ldots,P_n,Q)$ is ({\em totally}) {\em
$\Phi$-preserved}, or $\Phi$ is ({\em totally}) {\em
$(P_1,\ldots,P_n,Q)$-preserving}, if for any $\overline{a}_1\in
P_1,\ldots,\overline{a}_n\in P_n$,
$$\Phi(\overline{a}_1,\ldots,\overline{a}_n,\mathcal{M})\subseteq
Q.$$ Here we also say on {\em universal} $\Phi$- and
$(P_1,\ldots,P_n,Q)$-preservation.

If
$\Phi(\overline{a}_1,\ldots,\overline{a}_n,\mathcal{M})\subseteq
Q$ for some $\overline{a}_1\in P_1,\ldots,\overline{a}_n\in P_n$,
then we say that $(P_1,\ldots,P_n,Q)$ is {\em existentially
$\Phi$-preserved}, or $\Phi$ is {\em existentially
$(P_1,\ldots,P_n,Q)$-preserving}.

If $\Phi(\overline{a}_1,\ldots,\overline{a}_n,\mathcal{M})\cap
Q\ne\emptyset$ for some $\overline{a}_1\in
P_1,\ldots,\overline{a}_n\in P_n$ 
then we say that $(P_1,\ldots,P_n,Q)$ is {\em $\exists$-partially
$\Phi$-preserved}, or $\Phi$ is {\em $\exists$-partially
$(P_1,\ldots,P_n,Q)$-preserving}. If this property holds for any
$\overline{a}_1\in P_1,\ldots,\overline{a}_n\in P_n$, we say that
$(P_1,\ldots,P_n,Q)$ is {\em $\forall$-partially
$\Phi$-preserved}, or $\Phi$ is {\em $\forall$-partially
$(P_1,\ldots,P_n,Q)$-preserving}.

We say that the tuple $(P_1,\ldots,P_n,Q)$ is {\em
$\exists$-partially $\Phi$-non-preserved}, or $\Phi$ is {\em
$\exists$-partially $(P_1,\ldots,P_n,Q)$-non-preserving}, if
$\Phi(\overline{a}_1,\ldots,\overline{a}_n,\mathcal{M})\cap
\overline{Q}\ne\emptyset$ for some $\overline{a}_1\in
P_1,\ldots,\overline{a}_n\in P_n$, where
$\overline{Q}=M^m\setminus Q$. If this property holds for any
$\overline{a}_1\in P_1,\ldots,\overline{a}_n\in P_n$, we say that
$(P_1,\ldots,P_n,Q)$ is {\em $\forall$-partially
$\Phi$-non-preserved}, or $\Phi$ is {\em $\forall$-partially
$(P_1,\ldots,P_n,Q)$-non-preserving}. In the latter case we also
say that the tuple $(P_1,\ldots,P_n,Q)$ is {\em totally
$\Phi$-non-preserved}, or $\Phi$ is {\em totally $(P_1,\ldots,$
$P_n,$ $Q)$-non-preserving}.

If
$\Phi(\overline{a}_1,\ldots,\overline{a}_n,\mathcal{M})\subseteq
\overline{Q}$ for some $\overline{a}_1\in
P_1,\ldots,\overline{a}_n\in P_n$, then we say that
$(P_1,\ldots,P_n,Q)$ is {\em existentially $\Phi$-disjoint}, or
$\Phi$ is {\em existentially $(P_1,\ldots,P_n,Q)$-disjointing}.

If
$\Phi(\overline{a}_1,\ldots,\overline{a}_n,\mathcal{M})\subseteq
\overline{Q}$ for any $\overline{a}_1\in
P_1,\ldots,\overline{a}_n\in P_n$, then we say that
$(P_1,\ldots,P_n,Q)$ is {\em totally $\Phi$-disjoint} or {\em
universally $\Phi$-disjoint}, or $\Phi$ is {\em totally
$(P_1,\ldots,P_n,Q)$-disjointing}, or {\em universally
$(P_1,\ldots,P_n,Q)$-disjointing}.

If $\Phi$ is a singleton $\{\varphi\}$ then
totally/existentially/partially $\Phi$-(non-)preserved/disjoint
tuples are called {\em totally/existentially/partially
$\varphi$-{\rm (}non-{\rm )}preserved/disjoint}, respectively, and
$\varphi$ is {\em totally/existentially/partially
$(P_1,\ldots,P_n,$ $Q)$-{\rm (}non-{\rm )}preserving/disjointing}.

If $P_1=\ldots =P_n=Q$ then $(P_1,\ldots,P_n,Q)$-preserving type
$\Phi$ is called {\em $(P_1,\ldots,P_n,Q)$-idempotent} and
$(P_1,\ldots,P_n,Q)$ is {\em $\Phi$-idempotent}.

If $\Phi=\{\varphi\}$ then we replace $\Phi$ by $\varphi$ in the
definitions above.

\medskip
By the definition we have the following properties:

\begin{proposition}\label{prop_exs_total}
$1.$ If a type $\Phi$ is totally
$(P_1,\ldots,P_n,Q)$-preserving/disjointing and
$P_1\times\ldots\times P_n\ne\emptyset$ then $\Phi$ is
existentially $(P_1,\ldots,P_n,Q)$-preserving/disjointing.

$2.$ If a type $\Phi$ is $\forall$-partially
$(P_1,\ldots,P_n,Q)$-{\rm (}non-{\rm )}preserving and
$P_1\times\ldots\times P_n\ne\emptyset$ then $\Phi$ is
$\exists$-partially $(P_1,\ldots,P_n,Q)$-{\rm (}non-{\rm
)}preserving.
\end{proposition}

The following example shows that the converse for the items of
Proposition \ref{prop_exs_total} does not hold in general.

\begin{example}\label{ex_exist_non-tot}\rm
1. Let $P_1,P_2,Q$ be nonempty unary predicates such that $P_1$
and $Q$ are connected by a bijection $f_0$, and a function
$f\supset f_0$ maps $P_2$ into the complement $\overline{Q}$ of
$Q$. Then the formula $f(x)\approx y$ is totally
$(P_1,Q)$-preserving and totally $(P_2,Q)$-disjointing while this
formula is existentially $(P_1\cup P_2,Q)$-preserving and not
totally $(P_1\cup P_2,Q)$-preserving.

2. Let the previous example be extended by a constant
$c\in\overline{Q}$ and the relation $f_0$ be extended by the set
$\{\langle a,c\rangle\mid a\in P_1\}$. Then this extended relation
$R(x,y)$ witnesses the $\forall$-partial $(P_1,Q)$-preservation.
Extending $P_1$ by new element $c'$ and $R$ by the pair $\langle
c',c\rangle$ we obtain the witness of $\exists$-partial
$(P_1\cup\{c'\},Q)$-preservation which is not $\forall$-partial
$(P_1\cup\{c'\},Q)$-preserving. Similarly there is a
$\exists$-partial $(P_1\cup\{c'\},Q)$-non-preservation which is
not $\forall$-partial $(P_1\cup\{c'\},Q)$-non-preserving.
\end{example}

\begin{proposition}\label{prop_exs_total2}
For any type $\Phi$ and definable or non-definable relations
$P_1,\ldots,P_n,Q$ in a structure $\mathcal{M}$ the following
conditions are equivalent:

$1)$ $\Phi$ is totally/existentially
$(P_1,\ldots,P_n,Q)$-preserving;

$2)$ $\Phi$ is totally/existentially
$(P_1,\ldots,P_n,\overline{Q})$-disjointing.

\end{proposition}

Proof follows by the definition.

\begin{proposition}\label{prop_Cart_prod} For any
$(P_1,\ldots,P_n,Q)$ and $\Phi$, $(P_1,\ldots,P_n,Q)$ is
totally/existentially/par\-ti\-al\-ly $\Phi$-{\rm (}non-{\rm
)}preserved/disjoint iff $(P_1\times\ldots\times P_n,Q)$ is
totally/existentially/partially $\Phi$-{\rm (}non-{\rm
)}preserved/disjoint, where
$(\overline{x}_1,\ldots,\overline{x}_n,\overline{y})$ in $\Phi$ is
replaced by
$(\overline{x}_1\hat{\,}\ldots\hat{\,}\,\overline{x}_n,\overline{y})$.
\end{proposition}

Proof follows since the condition $\overline{a}_1\in
P_1,\ldots,\overline{a}_n\in P_n$ is equivalent to the condition
$(\overline{a}_1\hat{\,}\ldots\hat{\,}\,\overline{a}_n)\in
P_1\times\ldots\times P_n$. \endproof

\medskip
The following assertion shows that type-definable sets can be
described in terms of preservation.

\begin{proposition}\label{prop_type_def}
For any type
$\Phi=\Phi(\overline{x}_1,\ldots,\overline{x}_n,\overline{y})$ and
tuples $\overline{a}_1,\ldots,\overline{a}_n$ forming the set
$Q=\Phi(\overline{a}_1,\ldots,\overline{a}_n,\mathcal{M})$ the set
$Q$ is characterized by the following conditions: $\Phi$ is
$(\{\overline{a}_1\},\ldots,$ $\{\overline{a}_n\},Q)$-preserving
and $\Phi$ is not {\rm (}$\forall$-{\rm )} $\exists$-partially
$(\{\overline{a}_1\},\ldots,\{\overline{a}_n\},\overline{Q})$ for
the complement $\overline{Q}=M^{l(\overline{y})}\setminus Q$. In
particular, for a formula
$\varphi=\varphi(\overline{x}_1,\ldots,\overline{x}_n,\overline{y})$
and tuples $\overline{a}_1,\ldots,\overline{a}_n$ forming the set
$Q=\varphi(\overline{a}_1,\ldots,\overline{a}_n,\mathcal{M})$ the
set $Q$ is characterized by the following conditions: $\varphi$ is
$(\{\overline{a}_1\},\ldots,\{\overline{a}_n\},Q)$-preserving and
$\neg\varphi$ is
$(\{\overline{a}_1\},\ldots,\{\overline{a}_n\},\overline{Q})$-preserving.
\end{proposition}

Proof follows by the definition. \endproof

\begin{remark}\label{rem_repr}\rm
Using Proposition \ref{prop_type_def} we have the following
representation of a type-definable set
$Q=\Phi(\overline{a}_1,\ldots,\overline{a}_n,\mathcal{M})$:
$$Q=\{\overline{b}\mid\Phi(\overline{x}_1,\ldots,\overline{x}_n,\overline{y})
\mbox{ is }(\forall\mbox{-})\mbox{ }\exists\mbox{-partially }
(\{\overline{a}_1\},\ldots,\{\overline{a}_n\},\{\overline{b}\})\mbox{-preserving}\}.$$
\end{remark}

\medskip
{\bf Definition} \cite{rs2, rs3, rs4}. Let $T$ be a complete
theory, $\mathcal{M}\models T$. Consider types
$p(\overline{x}),q(\overline{y})\in S(\emptyset)$, realized in
$\mathcal{M}$, and all {\em
$(p,q)$-preserving}\index{Formula!$(p,q)$-preserving} formulae
$\varphi(\overline{x},\overline{y})$ of $T$, i.~e., formulae for
which there is $\overline{a}\in M$ such that $\models
p(\overline{a})$ and $\varphi(\overline{a},\overline{y})\vdash
q(\overline{y})$. For each such a formula
$\varphi(\overline{x},\overline{y})$, we define a relation
$R_{p,\varphi,q}\rightleftharpoons\{(\overline{a},\overline{b})\mid\mathcal{M}\models
p(\overline{a})\wedge\varphi(\overline{a},\overline{b})\}.$ If
$(\overline{a},\overline{b})\in R_{p,\varphi,q}$, then the pair
$(\overline{a},\overline{b})$ is called a {\em
$(p,\varphi,q)$-arc}.

\medskip
The definition above is naturally spread for types
$\Phi(\overline{x},\overline{y})$ replacing the formula
$\varphi(\overline{x},\overline{y})$.

\begin{proposition}\label{prop_pres_tuple_pair}
For any types $p(\overline{x}),q(\overline{y})\in S(\emptyset)$
and a formula $\varphi(\overline{x},\overline{y})$ the following
conditions are equivalent:

$1)$ the formula $\varphi$ is $(p,q)$-preserving;

$2)$ the pair $(p(\mathcal{M}),$ $q(\mathcal{M}))$ is totally
$\varphi$-preserved;

$3)$ the pair $(p(\mathcal{M}),$ $q(\mathcal{M}))$ is
existentially $\varphi$-preserved.
\end{proposition}

Proof. $1)\Leftrightarrow 3)$ and $2)\Rightarrow 3)$ follow by the
definition.

$3)\Rightarrow 2)$ holds since any realizations $\overline{a}$ and
$\overline{a}'$ of the complete type $p(\overline{x})$ are
connected by an automorphism $f$ of an elementary extension of
$\mathcal{M}$ implying that
$\varphi(\overline{a},\overline{y})\vdash q(\overline{y})$ iff
$\varphi(\overline{a}',\overline{y})\vdash q(\overline{y})$, i.e.
$(\overline{a},\overline{b})\in p(\mathcal{M})\times
q(\mathcal{M})$ is a $(p,\varphi,q)$-arc iff
$(\overline{a}',f(\overline{b}))$ is a $(p,\varphi,q)$-arc.
\endproof

\medskip
By Propositions \ref{prop_Cart_prod} and
\ref{prop_pres_tuple_pair} we observe a connection between
syntactic and semantic conditions of preservation.

\medskip
The following assertion allows to reduce the type preservation for
type-definable properties till a formulaic one.

\begin{proposition}\label{prop_rem22}
If some conjunction of formulae in a type $\Phi$ is
totally/existentially $(P_1,\ldots,$ $P_n,Q)$-preserving/disjoint
then $\Phi$ is totally/exis\-t\-entially
$(P_1,\ldots,P_n,Q)$-preserving/disjoint.
\end{proposition}

Proof. If some conjunction $\varphi_0$ of formulae in $\Phi$ is
totally/existentially $(P_1,\ldots,P_n,Q)$-preserving/disjoint
then $\Phi$ is totally/exis\-t\-en\-ti\-al\-ly
$(P_1,\ldots,P_n,Q)$-preserving/disjoint since the set of
solutions for $\Phi$ is a subset of the set of solutions for
$\varphi_0$.
\endproof

The following assertion shows that the converse of Proposition
\ref{prop_rem22} can fail for type-definable properties.

\begin{proposition}\label{prop_rem222}
Let $p(\overline{x})$ and $q(\overline{y})$ be nonisolated types
forced by their sets of formulae $\psi_i(\overline{x})$,
$\chi_i(\overline{y})$, respectively, $i\in \omega$, such that
$$\vdash(\psi_i(\overline{x})\to\psi_j(\overline{x}))\wedge\exists\overline{x}(\psi_j(\overline{x})\wedge\neg\psi_i(\overline{x})),$$
$$\vdash(\chi_i(\overline{x})\to\chi_j(\overline{x}))\wedge\exists\overline{x}(\chi_j(\overline{x})\wedge\neg\chi_i(\overline{x})),$$
for any $i<j\in\omega$; $\Phi(\overline{x},\overline{y})$ consist
of formulae $\phi_i(\overline{x},\overline{y})$ such that for any
realization $\overline{a}\models p(\overline{x})$ and any
$j\in\omega$,
$$\models\exists\overline{y}(\phi_i(\overline{x},\overline{y})\wedge\chi_j(\overline{y})\wedge\neg\chi_{j+1}(\overline{y}))$$
iff $j\geq y$. Then for any saturated model $\mathcal{M}$ of the
given theory, $\Phi$ is totally
$(p(\mathcal{M}),q(\mathcal{M}))$-preserving whereas each finite
conjunction $\theta$ of formulae in $\Phi$ is $\forall$-partially
$(p(\mathcal{M}),q(\mathcal{M}))$-preserving and
$\forall$-partially
$(p(\mathcal{M}),q(\mathcal{M}))$-non-preserving, in particular,
$\theta$ is both not totally
$(p(\mathcal{M}),q(\mathcal{M}))$-preserving.
\end{proposition}

\medskip
Similarly to Proposition \ref{prop_rem22} we have the following:

\begin{proposition}\label{prop_rem23}
If a type $\Phi$ is $\alpha$-partially $(P_1,\ldots,P_n,Q)$-{\rm
(}non-{\rm )}preserving, where $\alpha\in\{\forall,\exists\}$,
then any conjunction of formulae in a type $\Phi$ is
$\alpha$-partially $(P_1,\ldots,P_n,Q)$-{\rm (}non-{\rm
)}preserving.
\end{proposition}

Proof follows immediately by the definition. \endproof

\begin{proposition}\label{prop_rem222a}
If properties $P_1,\ldots,P_n,Q$ are type-definable in a saturated
structure then a type $\Phi$ is par\-ti\-al\-ly
$(P_1,\ldots,P_n,Q)$-preserving/disjoint iff some conjunction of
formulae in $\Phi$ is totally/existentially/par\-ti\-al\-ly
$(P_1,\ldots,P_n,Q)$-preserving/disjoint.
\end{proposition}

\medskip
We have the following basic properties for variations of
preservation:

\begin{proposition}\label{pr21} {\rm (Monotony)}
If $(P_1,\ldots,P_n,Q)$ is $\Phi$-preserved, $P_1\supseteq
P'_1,\ldots,P_n\supseteq P'_n$, $Q\subseteq Q'$,
$\Phi\subseteq\Phi'$ then $(P'_1,\ldots,P'_n,Q')$ is
$\Phi'$-preserved.
\end{proposition}

Proof follows immediately by the definition. \endproof

\medskip
For types $\Phi$ and $\Psi$ we denote by $\Phi\vee\Psi$ the type
$\{\varphi\vee\psi\mid\varphi\in\Phi,\psi\in\Psi\}$, and by
$\Phi\wedge\Psi$ the type $\Phi\cup\Psi$, if it is consistent.
Since $\Phi\vee\Psi$ (respectively, $\Phi\wedge\Psi$) has the set
of solutions represented as the union (intersection) of the
type-definable sets for $\Phi$ and $\Psi$, we obtain the
following:

\begin{proposition}\label{pr22} {\rm (Union)}
If $(P_1,\ldots,P_n,Q)$ is $\Phi$-preserved and
$(P_1,\ldots,P_n,Q')$ is $\Psi$-pre\-ser\-v\-ed, with
$Q,Q'\subseteq M^m$, then $\Phi\vee\Psi$ is $(P_1,\ldots,P_n,Q\cup
Q')$-preserving and $\Phi\wedge\Psi$, if it is consistent, is
$(P_1,\ldots,P_n,Q\cap Q')$-preserving.
\end{proposition}

Proposition \ref{pr22} immediately implies the following:

\begin{corollary}\label{cor23}
If there is a $(P_1,\ldots,P_n,Q)$-preserving type $\Phi$ then the
family $Z_\Phi(P_1,\ldots,P_n,$ $Q)$ of all sets of solutions, in
a given structure $\mathcal{M}$, for
$(P_1,\ldots,P_n,Q)$-preserving types, which are contained in
$\Phi$, forms a distributive lattice $\langle
Z_\Phi(P_1,\ldots,P_n,Q);\cup,\cap\rangle$ with the least element
$\Phi(\mathcal{M})$.
\end{corollary}

\begin{remark}\label{rem24}\rm The lattice $\langle
Z_\Phi(P_1,\ldots,P_n,Q);\cup,\cap\rangle$ can have or not have
the greatest element $\Psi(\mathcal{M})$ depending on the
existence of greatest union of type-definable subsets
$\Phi(\overline{a}_1,\ldots,$ $\overline{a}_n,\mathcal{M})$ of
$Q$, where $\overline{a}_1\in P_1,\ldots,\overline{a}_n\in P_n$.

Besides, if the condition for subsets of $\Phi$ is omitted then
the union of considered types can be inconsistent violating that
$\cap$ is an operation here.
\end{remark}

\begin{remark}\label{rem25}\rm
The notions and results above can be naturally spread to the
families $\Phi$ of formulae with unboundedly many free variables.
So we can assume that $\Phi$ has arbitrarily many free variables.
In particular, considering sets $\Phi$ of formulae whose variables
belong to a countable set, one can admit countably many free
variables.
\end{remark}

\section{Traces of types}

{\bf Definition.} (cf. \cite{Newel}) Let
$\Phi=\Phi(\overline{x}_1,\ldots,\overline{x}_n,\overline{y})$ be
a type with consistent
$\Phi(\overline{a}_1,\ldots,\overline{a}_n,\overline{y})$, where
$\overline{a}_1,\ldots,\overline{a}_n$ be tuples in a model
$\mathcal{M}$ of a given theory $T$. A {\em trace} of $\Phi$ with
respect to $(\overline{a}_1,\ldots,\overline{a}_n)$, or a
$\Phi$-{\em trace}, is a family $\{Q_i\subseteq
M^{l(\overline{y})}\mid i\in I\}$ such that
$\Phi(\overline{a}_1,\ldots,\overline{a}_n,\mathcal{M})\subseteq\bigcup\limits_{i\in
I}Q_i$ and $\Phi$ is ($\forall$-) $\exists$-partially
$(\{\overline{a}_1\},\ldots,\{\overline{a}_n\},Q_i)$-preserving
for each $i\in I$.

If the sets $Q_i$ are pairwise disjoint then the $\Phi$-trace
$\{Q_i\mid i\in I\}$ is {\em disjoint}, too.

The $\Phi$-trace $\{Q_i\mid i\in I\}$ is called $A$-({\em
type}-){\em definable} if each $Q_i$ is a (type-)definable set,
which are defined over $A$. We say on the (type-)definability of
the trace is it is $A$-({\em type}-){\em definable} for some $A$.

\medskip
Notice that any type
$\Phi(\overline{x}_1,\ldots,\overline{x}_n,\overline{y})$ has a
type-definable trace
$\{\Phi(\overline{a}_1,\ldots,\overline{a}_n,\mathcal{M})\}$ over
the set $\cup\overline{a}_1\cup\ldots\cup\cup\overline{a}_n$,
which is a singleton. Similarly each type $\Theta(\overline{y})$
with
$\Phi(\overline{a}_1,\ldots,\overline{a}_n,\mathcal{M})\vdash\Theta(\overline{y})$
produce a singleton-trace $\{\Theta(\mathcal{M})\}$ for $\Phi$. By
the definition all these traces are disjoint.

\begin{example}\label{ex_dis_trace}\rm
Each type
$\Psi=\Phi(\overline{a}_1,\ldots,\overline{a}_n,\overline{y})$ has
a definable trace
$\{\{\overline{b}\}\mid\,\,\models\Phi(\overline{a}_1,\ldots,\overline{a}_n,\overline{b})\}$.
Here the trace is
$\Phi(\overline{a}_1,\ldots,\overline{a}_n,\mathcal{M})$-definable,
and it is $\emptyset$-definable iff
$\Phi(\overline{a}_1,\ldots,\overline{a}_n,\mathcal{M})\subseteq{\rm
dcl}(\emptyset)$.
\end{example}

\medskip
For a type
$\Psi=\Phi(\overline{a}_1,\ldots,\overline{a}_n,\overline{y})$ we
denote by $[\Psi]$ the set $$\{p(\mathcal{M})\mid
p(\overline{y})\in
S^{l(\overline{y})}(\emptyset),\Phi(\overline{a}_1,\ldots,\overline{a}_n,\mathcal{M})\cap
p(\mathcal{M})\ne\emptyset\}.$$

\begin{proposition}\label{prop_tr_disj}
For any type
$\Psi=\Phi(\overline{a}_1,\ldots,\overline{a}_n,\overline{y})$ the
family $[\Psi]$ is a disjoint $\emptyset$-type-definable
$\Phi$-trace. It is definable iff each type $p(\overline{y})$ for
$[\Psi]$ is isolated.
\end{proposition}

Proof. If
$\overline{b}\in\Phi(\overline{a}_1,\ldots,\overline{a}_n,\mathcal{M})$
then $\overline{b}\in p(\mathcal{M})$, where $p(\overline{y})={\rm
tp}(\overline{b})$. Therefore $\Phi$ is ($\forall$-)
$\exists$-partially
$(\{\overline{a}_1\},\ldots,\{\overline{a}_n\},Q_i)$-preserving
and $\overline{b}\in\bigcup[\Psi]$. Hence $[\Psi]$ is a trace of
$\Phi$ with respect to $(\overline{a}_1,\ldots,\overline{a}_n)$.
Since the trace $[\Psi]$ collects sets of realizations of complete
types, which are pairwise inconsistent, it is disjoint. The trace
$[\Psi]$ is type-definable by the definition.

The latter characterization of the definability follows since the
non-empty set of realizations of a complete type $q(\overline{y})$
is definable iff $q(\overline{y})$ is isolated by a formula
$\varphi(\overline{y})$. Here
$q(\mathcal{M}=\varphi(\mathcal{M})$.
\endproof

\medskip
Proposition \ref{prop_tr_disj} immediately implies:

\begin{corollary}\label{cor_tr1}
If the model $\mathcal{M}$ is atomic then each $\Phi$-trace
$[\Psi]$ is definable.
\end{corollary}

Using Ryll-Nardzewski Theorem we have:

\begin{corollary}\label{cor_tr2}
If the model $\mathcal{M}$ of a countable language is saturated
then each $\Phi$-trace $[\Psi]$ is definable iff ${\rm
Th}(\mathcal{M})$ is $\omega$-categorical.
\end{corollary}

\begin{proposition}\label{prop_trace_singl}
For any type
$\Psi=\Phi(\overline{a}_1,\ldots,\overline{a}_n,\overline{y})$
with a saturated model $\mathcal{M}$ the trace $[\Psi]$ is finite
iff $\Phi(\overline{a}_1,\ldots,$ $\overline{a}_n,\overline{y})$
forces finitely many type in $S_{\overline{y}}(\emptyset)$.
\end{proposition}

Proof. If $[\Psi]$ is finite then the types for $[\Psi]$ are
divided by finitely many formulae
$\varphi_i=\varphi_i(\overline{y})$, $i<n$, into singletons
$[\Psi\cup\{\varphi_i\}]$. This singletons are consistent with
unique types $p_i(\overline{y})\in S(\emptyset)$, since
$\mathcal{M}$ is saturated. Hence $\Psi\cup\{\varphi_i\}$ forces
$p_i$ for each $i$ and $\Phi(\overline{a}_1,\ldots,$
$\overline{a}_n,\overline{y})$ forces
$\bigvee\limits_{i}p_i(\overline{y})$. The converse direction is
obvious. \endproof

\medskip
Proposition \ref{prop_trace_singl} immediately implies:

\begin{corollary}\label{cor_sing1} For any type
$\Psi=\Phi(\overline{a}_1,\ldots,\overline{a}_n,\overline{y})$
with a saturated model $\mathcal{M}$ the trace $[\Psi]$ is a
singleton iff $\Phi(\overline{a}_1,\ldots,$
$\overline{a}_n,\overline{y})$ forces a unique type in
$S_{\overline{y}}(\emptyset)$.
\end{corollary}

\begin{corollary}\label{cor_sing2}
If a type $\Phi(\overline{x}_1,\ldots,$
$\overline{x}_n,\overline{y})$ forces a $({\rm
tp}(\overline{a}_1\hat{\,}\ldots\hat{\,}\,\overline{a}_n),{\rm
tp}(\overline{b}))$-preserving formula
$\varphi(\overline{x}_1,\ldots,\overline{x}_n,\overline{y})$ for
some/any
$\overline{b}\in\Phi(\overline{a}_1,\ldots,\overline{a}_n,\mathcal{M})$
then $[\Phi(\overline{a}_1,\ldots,\overline{a}_n,\overline{y})]$
is a singleton.
\end{corollary}

\medskip
The following example shows that Proposition
\ref{prop_trace_singl} can fail if the given structure
$\mathcal{M}$ is not saturated.

\begin{example}\label{singl_fail}\rm
Let $\mathcal{M}$ be a structure in a signature of unary
predicates $P_i$, $i\in\omega$, $R_1$, $R_2$ with two nonempty
signature parts $R_1$, $R_2$ such that $R_1\cup R_2=M$, $R_1$ is
not divided by $P_i$, and the predicates $P_i$ form a substructure
on $R_2$ such that these $P_i$ are independent: all $R_2(x)\wedge
P^{\delta_1}_{i_1}(x)\wedge\ldots\wedge P^{\delta_n}_{i_n}(x)$ are
consistent, $i_1<\ldots<i_n<\omega$,
$\delta_1,\ldots,\delta_n\in\{0,1\}$, $n\in\omega$. Clearly the
formula $R_1(x)$ forces a complete type $p_1(x)$ whereas $R_2(x)$
belongs to continuum-many $1$-types $p_2^\Delta(x)$,
$\Delta\in\{0,1\}^\omega$, containing $R_2(x)$ and
$P^{\Delta(i)}_i(x)$, $i\in\omega$. Additionally we assume that
for some $\Delta_0\in\{0,1\}^\omega$ the complete type
$p_2^{\Delta_0}(x)$ is omitted in $\mathcal{M}$.

Now we consider the type $\Phi(y)$ formed by formulae
$R_1(y)\vee\varphi(y)$, $\varphi(y)\in p_2^{\Delta_0}(y)$. We have
$[\Phi(y)]=\{p_1(\mathcal{M})\}$ whereas $\Phi(y)$ does not force
unique type in $S_{\overline{y}}(\emptyset)$.

We also notice that $\Phi(y)$ does not force $({\rm
tp}(\emptyset),p_1)$-preserving formulae.
\end{example}

\begin{example}\label{ex_trace_acyclic}\rm
Let $\Gamma=\langle M;R\rangle$ be a graph with $R\ne\emptyset$.
Consider the type $\Phi=\{R(x,y)\}$ and $\Psi=\Phi(a,y)$ for an
element $a\in M$. We have $[\Psi]\ne\emptyset$ iff
$\Gamma\models\exists y R(a,y)$, and nonempty $[\Psi]$ is finite
iff $R(a,\Gamma)$ has finitely many $a$-orbits with respect to the
automorphism groups in elementary extensions of $\Gamma$. Clearly,
$[\Phi(a,y)]=[\Phi(b,y)]$ for any $b$ with ${\rm tp}(b)={\rm
tp}(a)$, but if ${\rm tp}(b)\ne{\rm tp}(a)$ then the equality
$[\Phi(a,y)]=[\Phi(b,y)]$ can hold or fail. Indeed, If
$R(a,\Gamma)\cap R(b,\Gamma)$ is a singleton $\{d\}$ then
$[\Phi(a,y)]$ and $[\Phi(b,y)]$ are singletons, too, with
$[\Phi(a,y)]=[\Phi(b,y)]=\{{\rm tp}(d)\}$. Respectively,
$[\Phi(a,y)]\ne[\Phi(b,y)]$ if, for instance, $R(a,\Gamma)$
consists of elements without loops and $R(b,\Gamma)$ has elements
with loops.
\end{example}

Now we spread the notion of trace $[\Psi]$ above for a type
$\Phi=\Phi(\overline{x}_1,\ldots,\overline{x}_n,\overline{y})$
with respect to a tuple $(P_1,\ldots,P_n)$, where $P_i\subseteq
M^{l(\overline{x}_i)}$, $i=1,\ldots,n$, in the following way.

We put
$[\Phi]_{(P_1,\ldots,P_n)}=\bigcup\{[\Phi(\overline{a}_1,\ldots,\overline{a}_n,\overline{y})]\mid\overline{a}_1\in
P_1,\ldots,\overline{a}_n\in P_n\}$. We write $[\Phi]$ instead of
$[\Phi]_{(P_1,\ldots,P_n)}$ if all $P_i$ cover all admissible sets
of tuples for $\Phi$, i.e. all tuples
$\overline{a}_1,\ldots,\overline{a}_n$ with consistent
$\Phi(\overline{a}_1,\ldots,\overline{a}_n,\overline{y})$ belong
to correspondent $P_1,\ldots,P_n$.

\medskip
By the definition, we have the following monotony property.

\begin{proposition}\label{prop_trace_mon}
For any type $\Phi$ and relations $P_i\subseteq P'_i$, $i\leq n$,
in a structure $\mathcal{M}$,
$$[\Phi]_{(P_1,\ldots,P_n)}\subseteq[\Phi]_{(P'_1,\ldots,P'_n)}.$$
\end{proposition}

\begin{lemma}\label{lemma_trace_type}
For any type
$\Phi(\overline{x}_1,\ldots,\overline{x}_n,\overline{y})$ and
tuples $\overline{a}_1,\ldots,\overline{a}_n$,
$\overline{b}_1,\ldots,\overline{b}_n$ with
$l(\overline{a}_i)=l(\overline{b}_i)=l(\overline{x}_i)$,
$i=1,\ldots,n$, if ${\rm
tp}(\overline{a}_1\hat{\,}\ldots\hat{\,}\overline{a}_n)={\rm
tp}(\overline{b}_1\hat{\,}\ldots\hat{\,}\overline{b}_n)$ then
$$[\Phi(\overline{a}_1,\ldots,\overline{a}_n,\overline{y})]=[\Phi(\overline{b}_1,\ldots,\overline{b}_n,\overline{y})].$$
\end{lemma}

Proof. Since  ${\rm
tp}(\overline{a}_1\hat{\,}\ldots\hat{\,}\overline{a}_n)={\rm
tp}(\overline{b}_1\hat{\,}\ldots\hat{\,}\overline{b}_n)$ there is
an automorphism $f$ of an elementary extension of given structure
such that $f(\overline{a}_i)=\overline{b}_i$, $i=1,\ldots,n$. This
automorphism moves
$\Phi(\overline{a}_1,\ldots,\overline{a}_n,\overline{y})$ onto
$\Phi(\overline{b}_1,\ldots,\overline{b}_n,\overline{y})$. Since
all automorphisms preserves types over $\emptyset$, we have
$[\Phi(\overline{a}_1,\ldots,\overline{a}_n,\overline{y})]=[\Phi(\overline{b}_1,\ldots,\overline{b}_n,\overline{y})].$
\endproof

\medskip
In view of Lemma \ref{lemma_trace_type} and the monotony property
we have:

\begin{corollary}\label{cor_trace_type}
For any type $\Phi=\Phi(\overline{x},\overline{y})$, a tuple
$\overline{a}$, its type $p(\overline{x})={\rm tp}(\overline{a})$,
and $\emptyset\ne P\subseteq p(\mathcal{M})$,
$[\Phi(\overline{a},\overline{y})]=[\Phi]_P$.
\end{corollary}

\begin{corollary}\label{cor_trace_type2}
For any type $\Phi=\Phi(\overline{x},\overline{y})$, a tuple
$\overline{a}$, its type $p(\overline{x})={\rm tp}(\overline{a})$,
a type $q(\overline{y})\in S(emptyset)$ and $\emptyset\ne
P\subseteq p(\mathcal{M})$, where $\mathcal{M}$ is saturated,
$\Phi$ is $(p,q)$-preserving iff
$[\Phi(\overline{a},\overline{y})]=[\Phi]_P=\{q(\overline{y})\}$.
\end{corollary}

{\bf Definition.} \cite{Kulp} A family $p_1,\ldots,p_n$ of
complete $1$-types over $A$ is {\em weakly orthogonal} over $A$ if
every $n$-tuple $\langle a_1,\ldots,a_n\rangle\in
p_1(\mathcal{M})\times\ldots\times p_n(\mathcal{M})$ satisfies the
same type over $A$. Here we omit $A$ if it is empty.

\medskip
In view of Corollary \ref{cor_trace_type} we have:

\begin{corollary}\label{cor_trace_type3}
For any type $\Phi=\Phi(x_1,\ldots,x_n,\overline{y})$ containing
$p_1(x_1)\cup\ldots\cup p_n(x_n)$ for weakly orthogonal family
$p_1(x_1),\ldots, p_n(x_n)$, and for any tuple $\langle
a_1,\ldots,a_n\rangle$,
$[\Phi(a_1,\ldots,a_n,\overline{y})]=[\Phi]_{(p_1(\mathcal{M}),\ldots,p_n(\mathcal{M}))}$.
\end{corollary}

For types
$\Phi_i=\Phi_i(\overline{x}_1,\ldots,\overline{x}_n,\overline{y}_i)$,
$i=1,\ldots,m$,
$\Psi=\Psi(\overline{y}_1,\ldots,\overline{y}_n,\overline{z})$ we
denote the type $S(\Phi_1,\ldots,\Phi_m,\Psi)$ consisting of all
formulae
\begin{equation}\label{eq_exist}
\exists
\overline{y}_1,\ldots,\overline{y}_m(\varphi_1(\overline{x}_1,\ldots,\overline{x}_n,\overline{y}_1)\wedge\ldots\varphi_m(\overline{x}_1,\ldots,\overline{x}_n,\overline{y}_m)\wedge
\psi(\overline{y}_1,\ldots,\overline{y}_m,\overline{z})),
\end{equation}
where
$\varphi_i(\overline{x}_1,\ldots,\overline{x}_n,\overline{y}_i)\in\Phi_i$,
$i=1,\ldots,m$,
$\psi(\overline{y}_1,\ldots,\overline{y}_m,\overline{z})\in \Psi$.

Notice that if the types $\Phi_1,\ldots,\Phi_m,\Psi$ define
operations $f_1,\ldots,f_m,g$, respectively, then the type
$S(\Phi_1,\ldots,\Phi_m,\Psi)$ defines the superposition
$g(f_1(\overline{x}_1,\ldots,\overline{x}_n),\ldots,f_m(\overline{x}_1,\ldots,\overline{x}_n),\ldots)$.
It correspond to the singleton trace with respect to the arguments
$\overline{a}_1,\ldots,\overline{a}_n$ for tuples
$\overline{x}_1,\ldots,\overline{x}_n$ of variables.

In general case the type $S(\Phi_1,\ldots,\Phi_m,\Psi)$ defines
the trace as a superposition of traces for
$\Phi_1,\ldots,\Phi_m,\Psi$, i.e. the superposition for operations
can be spread for type-definable relations:

\begin{theorem}\label{th_sup_tr}
For any types
$\Phi_i=\Phi_i(\overline{x}_1,\ldots,\overline{x}_n,\overline{y}_i)$,
$i=1,\ldots,m$,
$\Psi=\Psi(\overline{y}_1,\ldots,\overline{y}_n,\overline{z})$ and
tuples $\overline{a}_1,\ldots,\overline{a}_n$ for
$\overline{x}_1,\ldots,\overline{x}_n$, respectively, the trace
$[S(\Phi_1,\ldots,\Phi_m,\Psi)(\overline{a}_1,\ldots,\overline{a}_n,\overline{z})]$,
for a saturated structure $\mathcal{M}$, consists of all types
$p(\overline{z})\in S^{l(\overline{z})}(\emptyset)$ consistent
with $\Psi(\overline{b}_1,\ldots,\overline{b}_m,$ $\overline{z})$,
where ${\rm
tp}(\overline{b}_i)\in[\Phi_i(\overline{a}_1,\ldots,\overline{a}_n,\overline{y}_i)]$,
$i=1,\ldots,m$.
\end{theorem}

Proof. Let $p(\overline{z})\in S^{l(\overline{z})}(\emptyset)$ be
consistent with
$\Psi(\overline{b}_1,\ldots,\overline{b}_m,\overline{z})$, where
${\rm tp}(\overline{b}_i)\in[\Phi_i(\overline{a}_1,\ldots,$
$\overline{a}_n,\overline{y}_i)]$, $i=1,\ldots,m$, and
$\overline{c}\in\Psi(\overline{b}_1,\ldots,\overline{b}_n,\mathcal{M})\cap
p(\mathcal{M})$. Then the tuples
$\overline{b}_1,\ldots,\overline{b}_m$ witness the existences in
the formulae (\ref{eq_exist}). Hence
$p(\overline{z})\in[S(\Phi_1,\ldots,\Phi_m,\Psi)(\overline{a}_1,\ldots,\overline{a}_n,\overline{z})]$.

Conversely, if
$p(\overline{z})\in[S(\Phi_1,\ldots,\Phi_m,\Psi)(\overline{a}_1,\ldots,\overline{a}_n,\overline{z})]$
then by compactness there are tuples
$\overline{b}_1,\ldots,\overline{b}_m$ witnessing common
$\overline{y}_1,\ldots,\overline{y}_m$ in the formulae
(\ref{eq_exist}) and these tuples can be chosen in $\mathcal{M}$
since it is saturated. Hence $p(\overline{z})\in
S^{l(\overline{z})}(\emptyset)$ is consistent with
$\Psi(\overline{b}_1,\ldots,\overline{b}_m,\overline{z})$, where
${\rm
tp}(\overline{b}_i)\in[\Phi_i(\overline{a}_1,\ldots,\overline{a}_n,\overline{y}_i)]$,
$i=1,\ldots,m$. \endproof

\begin{remark}\label{rem_super}\rm
In view of Corollary \ref{cor_trace_type3} and Theorem
\ref{th_sup_tr} the trace
$$[S(\Phi_1,\ldots,\Phi_m,\Psi)(\overline{x}_1,\ldots,\overline{x}_n,\overline{z})]$$
is represented as the union of all traces
$[S(\Phi_1,\ldots,\Phi_m,\Psi)(\overline{a}_1,\ldots,\overline{a}_n,\overline{z})]$,
which are described in terms of superpositions (\ref{eq_exist}).
\end{remark}

\section{Preservations of properties by special formulae and
types}

\begin{remark}\label{rem31}\rm For a type $\Phi=\Phi(\overline{y})$ a tuple
$(P,\ldots,P_n,Q)$ is $\Phi$-preserved iff $Q\supseteq
\Phi(\mathcal{M})$, i.e. $Q$ contains that type-definable set.
Similarly, by the definition, in general case, $Q$ contains
type-definable sets
$\Phi(\overline{a}_1,\ldots,\overline{a}_n,\mathcal{M})$ whereas
$Q$ may be not type-definable.

If $\Phi$ contains sentences asserting that a binary relation $Q$
is (non, ir)reflexive, (anti)symmetric, transitive then the
correspondent type-definable set satisfies these properties. In
particular, types $\Phi$ allow to represent equivalence relations,
partial and linear orders.
\end{remark}

\begin{proposition}\label{pr32} If $\varphi$ is an atomic formula
$f(x_1,\ldots,x_n)\approx y$ then a tuple $(P_1,\ldots,$
$P_n,Q)=(P,\ldots,P,P)$ with $\emptyset\ne P\subseteq M$ is
$\varphi$-preserved, i.e. $\varphi$-idempotent, iff $P$ is the
universe of a subalgebra of a restriction of $\mathcal{M}$ till
the signature symbol $f$.
\end{proposition}

Proof follows by the definition since the restriction of
$\mathcal{M}$ till the signature symbol $f$ has a subalgebra with
the universe $P$ iff $P$ is closed under the operation $f$: for
any $a_1,\ldots,a_n\in P$, $f(a_1,\ldots,a_n)\in P$. Since the
value $f(a_1,\ldots,a_n)$ is unique, the condition
$f(a_1,\ldots,a_n)\in P$ means that both
$\varphi(x_1,\ldots,x_n,y)$ is
totally/$\forall$-partially/$\exists$-partially
$(\{a_1\},\ldots,\{a_n\},P)$-pre\-serv\-ing, i.e. the tuple
$(P,\ldots,P,P)$ is $\varphi$-idempotent. \endproof

\medskip
Proposition \ref{pr32} immediately implies:

\begin{corollary}\label{cor33} If $\Phi$ is the family of all atomic formula
$f(x_1,\ldots,x_n)\approx y$ for any functional signature symbol
$f$ of a structure $\mathcal{M}$ then a tuple $(P_1,\ldots,$
$P_n,Q)=(P,\ldots,P,P)$ with $\emptyset\ne P\subseteq M$ is
$\Phi$-preserved, i.e. $\Phi$-idempotent, iff $P$ is the universe
of a substructure of $\mathcal{M}$.
\end{corollary}

\begin{remark}\label{rem_many_sorted}\rm
Considering many-sorted structures
$\langle\mathcal{M}_i;\Sigma\rangle_{i\in I}$ these structures and
their many-sorted substructures can be expressed in terms of
preservations using the formulae $M_{i_1}(x_1)\wedge\ldots\wedge
M_{i_n}(x_n)\wedge f(x_1,\ldots,x_n)\approx y\wedge M_j(y)$ for
many-sorted operations $f$ with
$\delta_f=M_{i_1}\times\ldots\times M_{i_n}$ and $\rho_f\subseteq
M_j$. It means that in such a case each atomic formula
$f(x_1,\ldots,x_n)\approx y$ is
$(M_{i_1},\ldots,M_{i_n},M_j)$-preserving iff the functional
symbol $f\in\Sigma$ is interpreted by a function $f\mbox{: }
M_{i_1}\times\ldots\times M_{i_n}\to M_j$,
\end{remark}

\begin{remark}\label{rem_part_tot}\rm
Since the formulae $f(x_1,\ldots,x_n)\approx y$ have unique
solutions with respect to values $f(a_1,\ldots,a_n)$ in a
structure it does not matter to distinguish between partial and
total preservations for these formulae.
\end{remark}

\begin{proposition}\label{pr34} If
$\varphi(x^1_1,x^2_1;\ldots;x_n^1,x_n^2;y^1,y^2)$ is a formula
$$E(y^1,y^2)\wedge f(x^1_1,\ldots,x^1_n)\approx y^1\wedge
f(x^2_1,\ldots,x^2_n)\approx y^2$$ then a tuple $(P_1,\ldots,P_n,$
$Q)=(E,\ldots,E,E)$ with an equivalence relation $E\subseteq M^2$
is $\varphi$-preserved, i.e. $\varphi$-idempotent, iff $E$ is a
congruence relation of a restriction of $\mathcal{M}$ till the
signature symbol $f$.
\end{proposition}

Proof follows by the definition since the restriction of
$\mathcal{M}$ till the signature symbol $f$ is coordinated with
respect to the equivalence relation $E$, i.e. $E$ is a congruence
relation iff the conditions $E(a^1_1,a^2_1),\ldots,E(a^1_n,a^2_n)$
imply $E(f(a^1_1,\ldots,a^1_n),f(a^2_1,\ldots,a^2_n))$, which is
witnessed by the formula $\varphi$ with respect to the tuple
$(P_1,\ldots,P_n,$ $Q)=(E,\ldots,E,E)$.
\endproof

\section{Constructions of models of theories on a base of preservations of properties}

Let $T_0$ be a consistent theory. Following Henkin construction
\cite{Henkin, ErPa, Hodges} we extend the signature $\Sigma(T_0)$
by new constants and extend the theory $T_0$ till a complete
theory $T$ such that if $\exists y\varphi(y)\in T$ then
$\varphi(c)\in T$ for some constant $c$. A {\em canonical model}
$\mathcal{M}$ of $T$ \cite{Hodges} is represented by all infinite
equivalence classes $[t]=\{q\mid (q\approx t)\in T\}$, where terms
$t,q$ do not have free variables. These equivalence classes are
also represented by infinite $[c]=\{c'\mid (c'\approx c)\in T\}$,
where $c,c'$ are constant symbols, and it means that each constant
symbol $c'\in[c]$ is interpreted in $\mathcal{M}$ by the element
$[c]$. Here for any $n$-ary functional symbol $f\in\Sigma(T_0)$,
$f([c_1],\ldots,[c_n])=[c]$ iff $f(c_1,\ldots,c_n)\approx c\in T$,
and for any $n$-ary relational symbol $R\in\Sigma(T_0)$, $\models
R([c_1],\ldots,[c_n])$ iff $R(c_1,\ldots,c_n)\in T$.

Now we consider a characterization for $\mathcal{M}$ to be a model
of $T$ in terms of preservations of properties. Each formula
$\varphi(y)$ above is represented in the form
$\varphi(c_1,\ldots,c_n,y)$, where $c_1,\ldots,c_n$ are all
constant symbols in the formula $\varphi$. The condition
$\varphi(c_1,\ldots,c_n,c)\in T$, for appropriate $c$, means that
$\varphi(x_1,\ldots,x_n,x)$ is {\rm (}$\forall$-{\rm )}
$\exists$-partially $(\{[c_1]\},\ldots,\{[c_n]\},M)$-preserving.

Thus the canonical model $\mathcal{M}$ for the theory $T$, and its
restriction $\mathcal{M}\upharpoonright\Sigma(T_0)\models T_0$ are
characterized by the following {\em preserving condition}: a
$\Sigma(T_0)$-formula $\varphi(x_1,\ldots,x_n,x)$ is {\rm
(}$\forall$-{\rm )} $\exists$-partially
$(\{[c_1]\},\ldots,\{[c_n]\},M)$-preserving whenever $\exists
x\varphi(c_1,\ldots,c_n,x)\in T$.

Thus we have the following:

\begin{theorem}\label{th_can_mod} For any expansion
$\mathcal{M}$ of a model $\mathcal{M}_0$ of a theory $T_0$ by
naming all elements by infinitely many constants the following
conditions are equivalent:

$(1)$ $\mathcal{M}$ satisfies the preserving condition;

$(2)$ $\mathcal{M}\upharpoonright\Sigma(T_0)$ satisfies the
preserving condition;

$(3)$ $\mathcal{M}$ is a canonical model of a completion of $T_0$.
\end{theorem}

Proof. $(1)\Leftrightarrow(2)$ holds by the definition since the
{\rm (})$\forall$-{\rm )} $\exists$-partial
$(\{[c_1]\},\ldots,\{[c_n]\},M)$-preservation is considered for
$\Sigma(T_0)$-formulae.

$(1)\Rightarrow(3)$ is satisfied since the preserving condition
$\varphi([c_1],\ldots,[c_n],\mathcal{M})\ne\emptyset$ for a
$\Sigma(T_0)$-formula $\varphi(x_1,\ldots,x_n,x)$ is transformed
syntactically to the formula $\exists
x\varphi(c_1,\ldots,c_n,x)\to\varphi(c_1,\ldots,c_n,c)$, where
$[c]\in\varphi([c_1],\ldots,[c_n],\mathcal{M})$.

$(3)\Rightarrow(2)$ follows by the arguments above for the
preserving condition.
\endproof

\begin{remark}\label{rem_gen}\rm
The arguments for the representation of Henkin construction in
terms of preservations can be naturally spread for generic
constructions \cite{rs3}, both semantic and syntactic, with
confirmations of consistent formulae $\exists x\varphi(x)$ by
appropriate elements. These possibilities allow to realize links
of tuples with respect to arbitrary admissible diagrams and to
collect generic structures with given consistent lists of
properties.
\end{remark}

\section{Preservations of properties and Tarski-Vaught test}

Recall \cite{TV, Hodges} that a substructure $\mathcal{N}=\langle
N;\Sigma\rangle$ of a structure $\mathcal{M}=\langle
M;\Sigma\rangle$ is called an \emph{elementary substructure}
(denoted by $\mathcal{N}\preccurlyeq\mathcal{M}$), if for any
formula $\varphi(x_1,\ldots,x_n)$ of the signature $\Sigma$ and
any elements $a_1,\ldots$, $a_n\in N$ the condition
$\mathcal{N}\models\varphi(a_1,\ldots,a_n)$ is equivalent to the
condition $\mathcal{M}\models\varphi(a_1,\ldots,a_n)$. Here the
structure $\mathcal{M}$ is called the \emph{elementary extension}
of $\mathcal{N}$. If $N\ne M$, we write
$\mathcal{N}\prec\mathcal{M}$ instead of
$\mathcal{N}\preccurlyeq\mathcal{M}$. If
$\mathcal{N}\subseteq\mathcal{M}$ and the condition
$\mathcal{N}\preccurlyeq\mathcal{M}$
($\mathcal{N}\prec\mathcal{M}$) fails, we write
$\mathcal{N}\not\preccurlyeq\mathcal{M}$ (respectively,
$\mathcal{N}\not\prec\mathcal{M})$.

\begin{theorem}\label{TVT} {\rm (Tarski-Vaught Test) \cite{TV, Hodges}} Let $\mathcal{N}$
 be a substructure of a structure $\mathcal{M}$ in a signature $\Sigma$. Then the following conditions are equivalent:

{\rm (1)} $\mathcal{N}$ is an elementary substructure of
$\mathcal{M}$;

{\rm (2)} for any formula $\varphi(x_1,\ldots,x_n,y)$ of the
signature $\Sigma$ and any elements $a_1,\ldots,a_n\in N$ if
$\mathcal{M}\models\exists y\,\varphi(a_1, \ldots, a_n,y)$ then
there is an element $b\in N$ such that
$\mathcal{N}\models\varphi(a_1,\ldots,a_n,b)$.
\end{theorem}

\begin{corollary}\label{cor42} Let $\mathcal{N}$
be a substructure of a structure $\mathcal{M}$ in a signature
$\Sigma$. Then the following conditions are equivalent:

{\rm (1)} $\mathcal{N}$ is an elementary substructure of
$\mathcal{M}$;

{\rm (2)} for any formula $\varphi(x_1,\ldots,x_n,y)$ of the
signature $\Sigma$ and any elements $a_1,\ldots,a_n\in N$ if
$\mathcal{M}\models\exists y\,\varphi(a_1, \ldots, a_n,y)$ then
$\varphi(x_1, \ldots, x_n,y)$ is {\rm (}$\forall$-{\rm )}
$\exists$-partially $(\{a_1\},\ldots,\{a_n\},N)$-preserving;

{\rm (3)} for any formula $\varphi(x_1,\ldots,x_n,y)$ of the
signature $\Sigma$ and any elements $a_1,\ldots,a_n\in N$ either
$\varphi(x_1, \ldots, x_n,y)$ is {\rm (}$\forall$-{\rm )}
$\exists$-partially $(\{a_1\},\ldots,\{a_n\},M)$-disjoint or
$\varphi(x_1, \ldots, x_n,y)$ is {\rm (}$\forall$-{\rm )}
$\exists$-partially $(\{a_1\},\ldots,\{a_n\},N)$-preserving;

{\rm (4)} for any finite type $\Phi(x_1,\ldots,x_n,y)$ of the
signature $\Sigma$ and any elements $a_1,\ldots,a_n\in N$ either
$\Phi(x_1, \ldots, x_n,y)$ is {\rm (}$\forall$-{\rm )}
$\exists$-partially $(\{a_1\},\ldots,\{a_n\},M)$-disjoint or
$\Phi(x_1, \ldots, x_n,y)$ is {\rm (}$\forall$-{\rm )}
$\exists$-partially $(\{a_1\},\ldots,\{a_n\},N)$-preserving.
\end{corollary}

Proof. Since $\{a_1\},\ldots,\{a_n\}$ are singletons, the
$\forall$-conditions are equivalent to the $\exists$-conditions.

$(1)\Leftrightarrow(2)$ follows by Theorem \ref{TVT} since the
existence of $b\in N$ with
$\mathcal{N}\models\varphi(a_1,\ldots,a_n,b)$ means that
$\varphi(x_1, \ldots, x_n,y)$ is $\exists$-partially
$(\{a_1\},\ldots,\{a_n\},N)$-preserving.

$(2)\Leftrightarrow(3)$ holds since the condition
$\mathcal{M}\models\neg\exists y\,\varphi(a_1, \ldots, a_n,y)$
means that $\varphi(x_1, \ldots, x_n,y)$ is $\exists$-partially
$(\{a_1\},\ldots,\{a_n\},M)$-disjoint.

$(3)\Leftrightarrow(4)$ is satisfied since the type $\Phi(x_1,
\ldots, x_n,y)$ is $\exists$-partially
$(\{a_1\},\ldots,\{a_n\},M)$-disjoint iff the formula
$\bigwedge\Phi(x_1, \ldots, x_n,y)$ is $\exists$-partially
$(\{a_1\},\ldots,\{a_n\},M)$-disjoint, and the type $\Phi(x_1,
\ldots, x_n,y)$ is $\exists$-partially
$(\{a_1\},\ldots,\{a_n\},N)$-preserving iff the formula
$\bigwedge\Phi(x_1, \ldots,$ $x_n,y)$ is $\exists$-partially
$(\{a_1\},\ldots,\{a_n\},N)$-preserving. \endproof

\begin{remark}\label{rem43}\rm In Corollary \ref{cor42} the
condition of finiteness of the type $\Phi$ is essential since
infinite types can be both realized in a structure and omitted in
its elementary substructure. At the same time in the condition (4)
in Corollary \ref{cor42} the type $\Phi$ can be replaced by
infinite one if $\mathcal{N}$ is saturated.
\end{remark}

\section{Multipartite and related graphs with preservation properties}

Below we use the standard graph-theoretic terminology \cite{EmMe,
SO1}.

\medskip
{\bf Definition.} \cite{EmMe} For a cardinality $\varkappa$, a
{\em $\varkappa$-partite} graph is a graph whose vertices are (or
can be) partitioned into $\varkappa$ disjoint independent sets,
i.e. sets without arcs connecting elements inside these sets.
Equivalently, it is a graph that can be colored with $\varkappa$
colors, so that no two endpoints of an arc have the same color.
When  $\varkappa = 2$ these are the {\em bipartite} graphs, when
$\varkappa = 3$ they are called the {\em tripartite} graphs, etc.

\begin{proposition}\label{th_mp} Let $\Gamma=\langle M;R\rangle$ be a
graph. The following conditions are equivalent:

$(1)$ $\Gamma$ is $\varkappa$-partite;

$(2)$ $M$ is divided into disjoint parts $P_i$, $i<\varkappa$,
such that the formula $R(x,y)$ is
$(P_i,\overline{P_i})$-preserving for any $i<\varkappa$.
\end{proposition}

Proof. By the definition, if $\Gamma$ is a $\varkappa$-partite
graph, the universe $M$ is divided into $\varkappa$ disjoint parts
$P_i$ such that all arcs $(a,b)\in R$ with $a\in P_i$ has $b\in
\overline{P_i}$. Therefore $R(x,y)$ is
$(P_i,\overline{P_i})$-preserving for any $i<\varkappa$.
Conversely, if $R(x,y)$ is $(P_i,\overline{P_i})$-preserving for
any $i<\varkappa$ then all parts $P_i$ are independent implying
that $\Gamma$ is $\varkappa$-partite. \endproof

\medskip
Proposition \ref{th_mp} immediately implies:

\begin{corollary}\label{cor_bp} Let $\Gamma=\langle M;R\rangle$ be a
graph. The following conditions are equivalent:

$(1)$ $\Gamma$ is bipartite;

$(2)$ there is $P\subseteq M$ such that the formula $R(x,y)$ is
$(P,\overline{P})$-preserving and $(\overline{P},P)$-preserving.
\end{corollary}

\begin{corollary}\label{cor_tp} Let $\Gamma=\langle M;R\rangle$ be a
graph. The following conditions are equivalent:

$(1)$ $\Gamma$ is tripartite;

$(2)$ there are disjoint $P_0,P_1,P_2\subseteq M$ such that
$M=P_0\cup P_1\cup P_2$ and the formula $R(x,y)$ is
$(P_i,\overline{P_i})$-preserving for each $i<3$.
\end{corollary}

\begin{proposition}\label{prop_gr1}
Let $\Gamma=\langle M;R\rangle$ be a graph. The following
conditions are equivalent:

$(1)$ $R=\emptyset$ {\rm (}respectively, $R=M^2${\rm )};

$(2)$ the formula $R(x,y)$ {\rm (}$\neg R(x,y)${\rm )} is
$(M,\emptyset)$-preserving;

$(3)$ the formula $R(x,y)$ {\rm (}$\neg R(x,y)${\rm )} is
$(M,M)$-disjoint.
\end{proposition}

Proof. Since $\overline{M^2}=\emptyset$ is suffices to consider
the case $R=\emptyset$. It means that $R(a,\Gamma)=\emptyset$ for
any $a\in M$, which is equivalent to
$R(a,\Gamma)\subseteq\emptyset$ for any $a\in M$, i.e. the
$(M,\emptyset)$-preservation of $R(x,y)$, and to $R(a,\Gamma)\cap
M=\emptyset$ for any $a\in M$, i.e. the $(M,M)$-disjointness of
$R(x,y)$. \endproof

\section{Preservations of properties in graphs via distances}

In this section we consider an illustration for preservations of
properties using formulae describing distances in graphs.

Let $\Gamma=\langle M;R\rangle$ be a graph, $\varphi_n(x,y)$ be a
formula saying that there exists s $(x,y)$-path of length $n$,
$n\in\omega$. Taking the formula
$\psi_n(x,y):=\bigvee\limits_{m\leq n}\varphi_m(x,y)$, saying that
there exists s $(x,y)$-path of length $\leq n$, we observe the
following:

\begin{lemma}\label{lem_path} For any vertex $a\in M$ and a
property $P\subseteq M$ the formula $\psi_n(x,y)$ is partially
{\rm (}totally{\rm )} $(\{a\},P)$-preserving iff some {\rm
(}any{\rm )} vertex in $\Gamma$, which is achieved from $a$ by a
path of length $\leq n$, belongs to $P$.
\end{lemma}

Proof. If $\psi_n(x,y)$ is partially {\rm (}totally{\rm )}
$(\{a\},P)$-preserving then for some (any) $b\in M$ with
$\models\psi_n(a,b)$ we have $b\in P$. Therefore some {\rm
(}any{\rm )} vertex in $\Gamma$, which is achieved from $a$ by a
path of length $\leq n$, belongs to $P$.

Conversely, let some {\rm (}any{\rm )} vertex in $\Gamma$, which
is achieved from $a$ by a path of length $\leq n$, belongs to $P$.
Then we have $\psi_n(a,\Gamma)\cap P\ne\emptyset$
($\psi_n(a,\Gamma)\subseteq P$) confirming that formula
$\psi_n(x,y)$ is partially {\rm (}totally{\rm )}
$(\{a\},P)$-preserving. \endproof

\medskip
For a graph $\Gamma=\langle M;R\rangle$ and an element $a\in M$
the set $\{b\mid\Gamma\models R(a,b)\vee R(b,a)\}$ is called the
{\em $R$-neighbourhood} of $a$, of the length $1$. Replacing $R$
by $R^n$ we obtain the {\em $R$-neighbourhood} of $a$, of the
length $n$. Similarly we have the notion of
$\varphi$-neighbourhood for an arbitrary binary formula
$\varphi(x,y)$ in the signature $\langle R\rangle$.

\begin{theorem}\label{th_comp}
Let $\Gamma=\langle M;R\rangle$ be a graph, $P\subseteq M$ be a
nonempty property. Then the following conditions are equivalent:

$(1)$ $P$ is totally preserved under the formula $R(x,y)\vee
R(y,x)$;

$(2)$ $P$ is totally preserved under the formula
$\psi_n(x,y)\vee\psi_n(y,x)$ with some/any $n\geq 1$;

$(3)$ $P$ is a union of connected components.
\end{theorem}

Proof. $(1)\Rightarrow(2)$. Let $a\in P$. By the conjecture $P$ is
closed under $R$-neighbourhoods of length $1$, in particular, $P$
contains all elements in $\Gamma$ connected with $a$ by $R$-edges.
Continuing the process $n$ times we observe that $P$ contains all
elements having the $R$-distance at most $n$. Therefore $P$ is
totally preserved under the formula $\psi_n(x,y)\vee\psi_n(y,x)$
with some/any $n\geq 1$.

$(2)\Rightarrow(1)$ is obvious since $\psi_n$-neighbourhoods of
elements contain their $\psi_1$-neighbourhoods which are equal to
$R$-neighbourhoods.

$(2)\Rightarrow(3)$ follows by Lemma \ref{lem_path} since by the
conjecture for arbitrary element $a\in P$ all elements connected
with $a$ by $R$-paths belong to $P$, too.

$(3)\Rightarrow(1)$. Let $a$ be an arbitrary element in $P$. Since
$P$ contains the connected component $C(a)$ including $a$, then
$P$ contains the $R$-neighbourhood of $a$ of length $1$. Thus $P$
is totally preserved under the formula $R(x,y)\vee R(y,x)$.
\endproof

\medskip
Since connected components are minimal sets among their unions, we
have the following:

\begin{corollary}\label{cor_comp}
Connected components of a graph $\Gamma=\langle M;R\rangle$  are
exactly minimal nonempty subsets of its universe $M$, which are
totally preserved under the formula $R(x,y)\vee R(y,x)$.
\end{corollary}

\begin{remark}\label{rem_comp}\rm Replacing the formulae $R(x,y)\vee R(y,x)$ and
$\psi_n(x,y)\vee\psi_n(y,x)$ by $R(x,y)\wedge R(y,x)$ and
$\psi_n(x,y)\wedge\psi_n(y,x)$ we obtain the possibility of mutual
achievements of vertices in $\Gamma$ but it is possibly not
sufficient with respect to strong components since mutual
connections can have distinct lengths of counter paths.

By these arguments both connected components and strong components
may be not definable but it is enough if diameters of components
are equal to natural numbers.
\end{remark}

The following example shows that preservations by the formulae
$\varphi_n(x,y)$ can reconstruct initial graphs after replacements
of edges by mutual paths of the fixed length $n$.

\begin{example}\label{ex_comp}\rm
Let $\Gamma=\langle M;R\rangle$ be an arbitrary undirected graph,
$\Gamma'\langle M';R'\rangle$ is obtained from $\Gamma$ by adding
$(a,b)$-paths of a length $n\geq 2$ with new intermediate elements
of degree $2$ for each $(a,b)\in R$ with removals all edges in
$R$, i.e. with $R'\cap R=\emptyset$, such that these new
intermediate elements are pairwise disjoint for distinct
$(a,b)$-paths. Thus $M'\supseteq M$ and $M'\setminus M$ consists
of new intermediate elements. We have totally
$\varphi_n(x,y)$-preserved set $M$ and a family of totally
$\varphi_n(x,y)$-preserved subsets $M_1,\ldots,M_{n-1}$ of
$M'\supseteq M$ generated by intermediate elements in fixed
distances $m<n$ from elements in $M$. The reconstruction of
connected components in $\Gamma$ means the choice of elements in
these connected components and the step-by-step closure under
solutions of the formulae $\varphi_n(b,y)$, where $b$ are these
chosen elements on the initial step and obtained elements on
previous steps for subsequent ones.
\end{example}

Recall that the {\em diameter} of a connected undirected graph
$\Gamma$ is the supremum of shortest $(a,b)$-paths for all
elements $a,b$. If the supremum is infinite it is denoted by
$\infty$. It is also denoted by $\infty$ if the graph is not
connected.

The following assertion characterizes the infinite diameter of a
connected undirected graph in terms of partial preservation by
formulae $\varphi_n(x,y)$.

\begin{proposition}\label{prop_diam}
For any connected undirected graph $\Gamma=\langle M;R\rangle$ the
following conditions are equivalent:

$(1)$ $\Gamma$ has the diameter $\infty$;

$(2)$ for any formula $\varphi_n(x,y)$, $n\in\omega$, there is a
vertex $a\in M$ such that $\varphi_n(x,y)$ is partially
$(\{a\},M)$-preserving.
\end{proposition}

Proof. $(1)\Rightarrow(2)$. Since $\Gamma$ is connected the
diameter $\infty$ means that for any $n\in\omega$ there are
vertices $a_n,b_n\in M$ such that the shortest $(a_n,b_n)$-path
has the length $n$. Therefore $\varphi_n(x,y)$ is partially
$(\{a_n\},M)$-preserving. Since $n$ is arbitrary, it confirms the
item $(2)$.

$(2)\Rightarrow(1)$ follows by the definition since the partial
$(\{a\},M)$-preservation by $\varphi_n(x,y)$ means that there is
$b\in M$ such that the shortest $(a,b)$-path has the length $n$.
As $n$ is unbounded then the diameter is equal to $\infty$.
\endproof

\medskip
\begin{remark}\label{rem2}\rm The arguments above can be naturally transformed for directed
graphs replacing arcs by edges and by replacements
$\varphi'_n(x,y)$ of $\varphi_n(x,y)$ using the formulae
$R(x',y')\vee R(y',x')$ instead of $R(x',y')$. Thus we obtain a
characterization of the infinite diameter for directed graphs
again in terms of partial preservation.

Besides the finite diameters are also characterized by the
property of non-partially $(\{a\},M)$-preservation by formulae
$\varphi_n(x,y)$ and $\varphi'_n(x,y)$ starting with some $n$.
\end{remark}

\section{Preservation properties for algebraic constructions and decompositions}

Using the assertions above one can describe series of
type-definable structures, in particular, classes of (ordered)
semigroups, groups, rings and fields, including spherically
ordered ones \cite{rs5}, bands of semigroups including rectangular
bands of groups \cite{CP}, graded algebras \cite{grad}, etc.,
their subalgebras and quotients.

There are many kinds of graded structures $\mathcal{M}$ united by
the following one: $\mathcal{M}$ contains a binary operation
$\cdot$ and it is divided into parts $X_i$, $i\in I$, such that
for any $i,j\in I$ there is $k\in I$ such that $X_i\cdot
X_j\subseteq X_k$. In the introduced terms of preservation it
means that the formulae $x_i\cdot x_j\approx x_k$ are
$(X_i,X_j,X_k)$-preserved.

These kinds of preservation admit a series of natural
generalizations both with respect to $n$-ary operations preserving
the formulae $f(x_{i_1},\ldots,x_{i_n})\approx x_j$ by tuples
$(X_{i_1},\ldots,X_{i_n}, X_j)$ of parts, and more general natural
$\Phi(\overline{x}_1,\ldots,\overline{x}_m,\overline{y})$-type
preservations, including total and partial ones.

These possibilities of preservations allow to decompose structures
$\mathcal{M}$ into families of substructures with universes $X_i$
assuming their idempotency preservation with respect to given
operations.

\end{document}